\documentclass[12pt,twoside]{article} 
\usepackage{amssymb,graphicx}
\usepackage[small,sc]{caption2} 

\newtheorem{theorem}{Theorem}[section]
\newtheorem{lemma}[theorem]{Lemma}
\newtheorem{proposition}[theorem]{Proposition}

\newtheorem{rmrk}[theorem]{Remark}

\makeatletter
\@addtoreset{equation}{section}
\makeatother

\newcommand{\fig}[3] {
\medskip\smallskip
\begin{figure}[htb]
  \centering
  \includegraphics[width=#2]{#1.eps}
  \begin{minipage}[t]{0.80\linewidth} 
    \caption{#3}
    \protect\label{#1}
  \end{minipage}
\end{figure}
\medskip
}

\newenvironment{remark}
{\begin{rmrk} \em}
{\end{rmrk}}

\newcommand{\bi} {billiard}
\newcommand{\fn} {function}
\newcommand{\me} {measure}
\newcommand{\tr} {trajector}
\newcommand{\erg} {ergodic}
\newcommand{\sy} {system}
\newcommand{\hyp} {hyperbolic}

\newcommand{\dsy} {dynamical system}

\newcommand{\R} {\mathbb{R}}

\newcommand{\Z} {\mathbb{Z}}

\newcommand{\qed} {\hfill {\small Q.E.D.} \par\medskip}
\newcommand{\skippar} {\par\medskip}
\newcommand{\ds} {\displaystyle}
\newcommand{\proof} {\noindent \textsc{Proof.} }
\newcommand{\proofof}[1] {\noindent \textsc{Proof of {#1}.} }
\newcommand{\article}[3] {\textsc{{#1}}, {\itshape {#2}}, {{#3}}.}
\newcommand{\book}[3] {\textsc{{#1}}, {\itshape {#2}}, {{#3}}.}
\newcommand{\vol} {\textbf}
\newcommand{\eps} {\varepsilon}
\newcommand{\rset}[2] {\left\{ #1 \: \left| \: #2 \right. \! \right\} }
\newcommand{\lset}[2] {\left\{ \left. \! #1 \: \right| \: #2 \right\} }

\renewcommand{\iff} {if and only if\ }

\newcommand{\ta} {\Omega}            % billiard table
\newcommand{\bo} {\Gamma}            % border of billiard table
\newcommand{\ps} {\mathcal{M}}       % phase space
\newcommand{\ma} {\mathcal{F}}       % phase space
\renewcommand{\a} {\alpha}           % angle
\newcommand{\si} {\mathcal{S}}       % singularity set

\begin{document}

\title{\textbf{Hyperbolic billiards with nearly flat \\ 
focusing boundaries. I}}

\author{
\textsc{Luca Bussolari} 
\thanks{
Department of Mathematical Sciences,
Stevens Institute of Technology, 
Hoboken, NJ 07030, U.S.A.} $^\ddagger$
\qquad
\textsc{Marco Lenci} 
$^*$\thanks{
Dipartimento di Matematica,
Universit\`a di Bologna, 
P.zza di Porta S.~Donato 5,
40126 Bologna, ITALY}
\thanks{E-mails: \texttt{lbussola@math.stevens.edu},
\texttt{lenci@dm.unibo.it}} }

\date{December 2007}

\maketitle

\begin{abstract}
  The standard Wojtkowski--Markarian--Donnay--Bunimovich technique for
  the hyperbolicity of focusing or mixed billiards in the plane
  requires the diameter of a billiard table to be of the same order as
  the largest ray of curvature along the focusing boundary. This is
  due to the physical principle that is used in the proofs, the
  so-called defocusing mechanism of geometrical optics. In this paper
  we construct examples of hyperbolic billiards with a focusing
  boundary component of arbitrarily small curvature whose diameter is
  bounded by a constant independent of that curvature. Our proof
  employs a nonstardard cone bundle that does not solely use the
  familiar dispersing and defocusing mechanisms. 

\end{abstract}

\section{Introduction}
\label{sec-intro}

Much has been written, in the scientific literature, about the \hyp
ity of \bi s in two dimensions.  So much that general principles have
even been devised for the `design of \bi s with nonvanishing Lyapunov
exponents'. The expression is taken from the title of the 1986 seminal
paper by Wojtkowski \cite{w2}, in which he beautifully links the
question of exponential instability (i.e., positivity of a Lyapunov
exponent) to a few simple observations from geometrical optics. By
means of the powerful \emph{invariant cone technique} \cite{w1, k,
cm}, Wojtkowski gives sufficient conditions for a planar \bi\ to
have nonzero Lyapunov exponents, this implying a fuller range of \hyp\
properties via the general results of Katok and Strelcyn on Pesin's
theory for \dsy s with singularities \cite{ks}.

Wojtkowski's conditions are rather undemanding for \emph{dispersing}
and \emph{semidispersing} \bi s (i.e., \bi s in a domain $\Omega
\subset \R^2$, a.k.a.\ \emph{table}, whose boundary is the finite
union of smooth convex pieces, when seen from inside $\Omega$), and
much more restrictive for \emph{focusing}, \emph{semifocusing} and
\emph{mixed} \bi s (that is, cases when $\partial \Omega$ is made
up---completely or partially, respectively---by concave pieces). (Both
in the dispersing and in the focusing case, the prefix semi- means
that $\partial \Omega$ has some flat parts as well.) For the latter
type of \bi s, further work has been done by Markarian \cite{m1,m2},
Donnay \cite{d} and Bunimovich \cite{b3} (see \cite[Chap.~9]{cm} for
an overview of the subject and \cite{de} for an interesting
variation).

If we call \emph{boundary component} each smooth piece of $\partial
\Omega$, one of the conditions in \cite{w2} is that the inner
semiosculating disc at any given point of a focusing boundary
component must not intersect other components, or the semiosculating
discs relative to other focusing components (\cite{m1} has a similar
condition). This is required in order to implement the so-called
\emph{defocusing mechanism}, which can be loosely described like this:
One wants diverging beams of \tr ies to keep diverging after every
collision with the boundary. But at a focusing portion of the boundary
a diverging beam may be bounced back as a converging beam. A solution
around this problem is to let the converging beam travel untouched for
a sufficienly long time until the \tr ies focus among themselves and
then start to diverge again.

The defocusing mechanism is the closest extension of Sinai's original
idea of extracting \hyp ity from the expanding features of dispersing
boundaries \cite{s}. At least to our knowledge, it has remained
unsurpassed since Bunimovich introduced it in 1974 \cite{b1}, to
become very popular a few years later, when it was used to work out
the famous stadium \bi\ \cite{b2}.

Sticking too much to the standard principles, however, creates a
problem and somehow a paradox. The condition on the semiosculating
discs, and each of its later analogues, requires a table with focusing
components to have a diameter of the order of the largest radius of
curvature among the focusing points of the boundary.  To illustrate
how this may seem a paradox, consider the following example: Take a
unit square and replace three of its sides with circular arcs of
curvature $k_d \in (-\sqrt{2},0)$ having their endpoints in the
vertices of the square. In this paper we use the convention that the
curvature is positive at focusing points of the boundary and negative
at dispersing points, so the arcs are convex relative to the interior
of the square; the condition $|k_d| < \sqrt{2}$ ensures that each pair
of adjacent arcs intersects only at the common endpoint. The resulting
\bi\ is semidispersing, thus belongs to the standard class and is
well-know to be uniformly \hyp, Bernoulli, and so on \cite{cm}. Now
perturb the fourth side into a focusing circular arc of curvature $k_f
\ll 1$.  Now matter how small the perturbation, this new \bi\ will
never satisfy Wojtkowski's principle and is not currently known to be
\hyp, although it presumably is.

This may not sound too strange. After all, certain perturbations of
dispersing \bi s are known to possess elliptic islands \cite{rt,
tr}. But the paradox is that the smaller the perturbation, the less
adequate the standard technique; that is, the closer the \bi\ comes to
be dispersing, the worse the method applies which is supposed to
exploit the dispersing nature of the boundaries.  Up until $k_f=0$, at
which point everything suddenly, and abruptly, works again to the
fullest power of the theory of \hyp\ \bi s. 

\skippar

Here we address this problem and, although we cannot yet prove that
the perturbed square \bi\ is \hyp, we devise a couple of models that
make clear what the difficulties are in extending the current
methodology. These \bi s, which are modifications of the example just
discussed, are depicted in Figs.~\ref{fig-t1intro} and
\ref{fig-t2intro}. They are indeed two families of \bi s, as we are
interested in the case when the curvature of the focusing boundary
goes to zero. We define an invariant cone bundle that exploits the
fact that the focusing component is nearly flat, and thus almost
always acts as a semidispersing boundary.

\fig{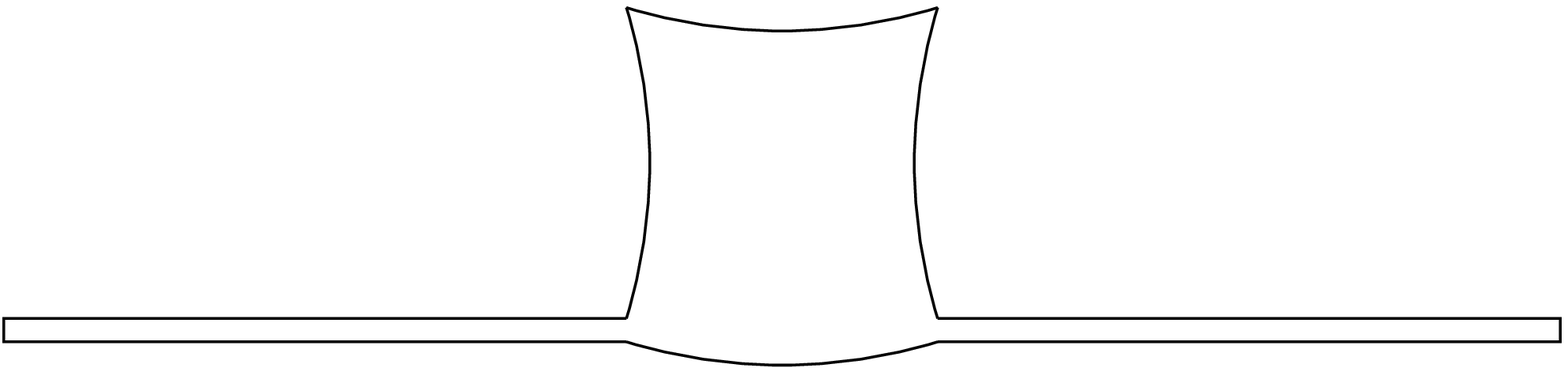} {12cm} {The main \bi\ table}

In any event, we are able to answer the following questions in the
affirmative:

\begin{enumerate}
\item Can one design a \bi\ whose \hyp ity is proved via a set of
  invariant cones that does not use exclusively the
  dispersing/defocusing mechanism for beams of \tr ies?

\item \label{pt2} Can one construct a family of \hyp\ \bi\ tables with
  a (nonvanishing) focusing component whose maximum curvature
  approaches zero, and such that the area of the table is bounded
  above?

\item \label{pt3} Can one require the diameter to be bounded above as
  well?

\item Are these \bi s \erg? (This will be proved in \cite{bl}.)
\end{enumerate}

\fig{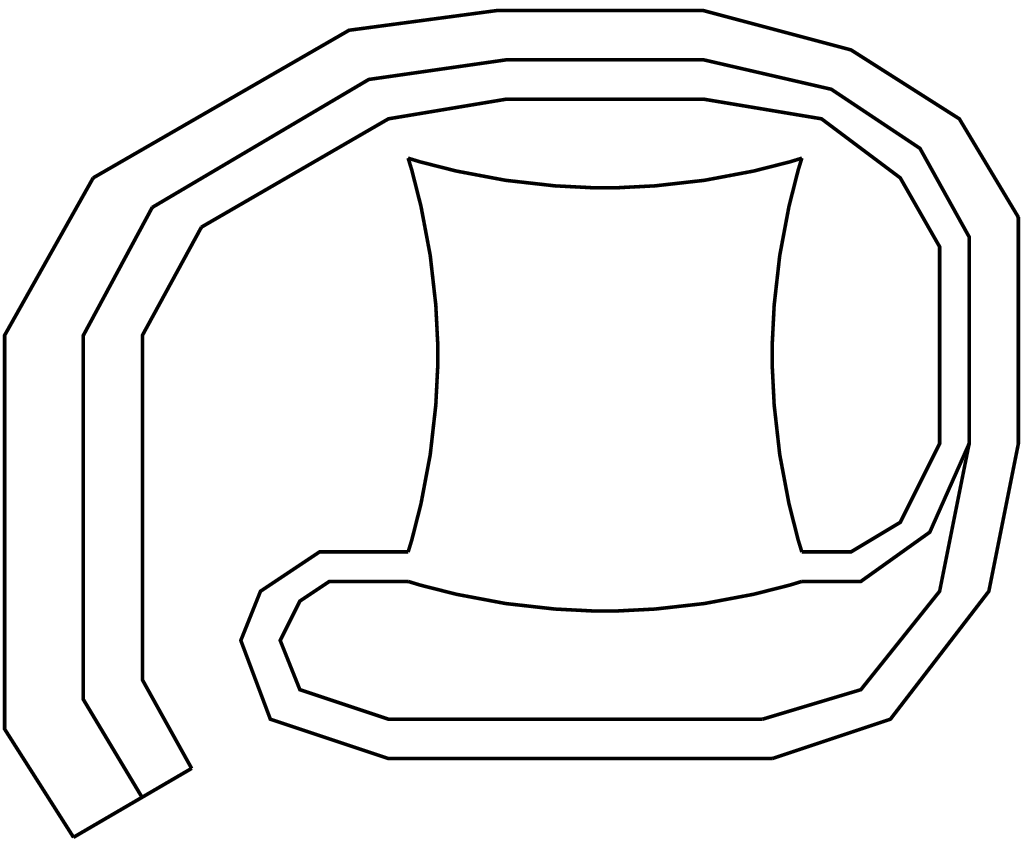} {6.2cm} {A modification of the main \bi\ table}

Points \ref{pt2} and \ref{pt3} show, independently of the method
utilized, that one can go beyond the apparent implication `almost flat
focusing boundaries imply very large tables'.

\skippar

This is the plan of the paper: In Section \ref{sec-prel} we review the
basic definitions of \bi\ dynamics. In Section \ref{sec-cones} we
present and adapt Wojtkowsky's theory of invariant cones derived from
geometrical optics. In Section \ref{sec-hyp} we define the first of
our models and choose suitable cones to prove its \hyp ity. In Section
\ref{sec-conf} we show that the \bi\ introduced before can be chosen
with a bounded area, and finally we present a second model which has
a bounded diameter as well.

\bigskip
\noindent
\textbf{Acknowledgments.} \ We would like to thank Gianluigi Del Magno
for an instructive discussion on the subject. M.L.\ acknowledges
partial support from NSF Grant DMS-0405439.

\section{Preliminaries}
\label{sec-prel}

A planar \bi\ is the \dsy\ generated by the flow of a point particle
that moves inertially inside a closed region $\ta \subset \R^2$ and
collides elastically at the boundary; the latter is assumed to have an
infinite mass. This implies that the \tr y of the particle, near the
collision point, verifies the well-known \emph{Fresnel law}: the angle
of incidence equals the angle of reflection.  The region $\ta$ is
called the \emph{\bi\ table}. We denote $\bo = \partial \ta$ and
assume that $\bo$ is piecewise smooth (at least $C^3$).

Let $(q(t),u(t))$ represent the position and the velocity of the
particle at time $t$. It is an easy consequence of the conservation of
energy that $\|u(t)\| =$ constant. Therefore, by a rescaling of time,
one can always reconduct to the situation where $\|u\| = 1$, which we
assume throughout the paper. The product $\ta \times S^1$ is the
natural phase space of the \bi\ flow, with a couple of extra
specifications: First, if $q \in \bo$ and $u$ points outwardly, then
$(q,u)$ is identified with $(q,u')$, where $u'$ is the outgoing (i.e.,
inward) velocity of a collision at $q$ with incoming velocity $u$.
Second, if $q$ is in a corner, the flow is not defined. The \bi\ flow
preserves the Lebesgue \me\ on $\ta \times S^1$, as it can be verified
directly or by applying the Liouville Theorem to this nonsmooth
Hamiltonian \sy.

Now let $\ps \subset \ta \times S^1$ be the set of all pairs $(q,u)$
with $q \in \bo$ and $u$ pointing inside the table. These pairs are
sometimes called \emph{line elements} \cite{s} and $\ps$ is evidently
a global cross section for the flow. The corresponding Poincar\'e map
$\ma: \ps \longrightarrow \ps$ is called the \emph{\bi\ map} and acts
as follows: if $q' \in \bo$ is the first collision point of the
flow-\tr y with initial conditions $(q,u)$, and $u'$ is the
postcollisional velocity there, then $\ma(q,u) = (q',u')$. $\ma (q,u)$
is undefined when $q'$ is a vertex of $\bo$, and is discontinuous at
tangential collisions, i.e., when $u'$ is tangent to $\bo$ in
$q'$. For the sake of simplicity, those latter line elements are
removed as well from the domain of $\ma$. The set of all removed
$(q,u)$ is denoted $\si_1$, or $\si_1^+$.

We identify $\ps$ with the rectangle $[0,L] \times [-\pi/2,\pi/2]$,
where $L$ is the perimeter of $\ta$: each $(q,u)$ is identified with
the pair $(s,\a)$, where $s$ is the arclength coordinate of $q$
(relative to a fixed choice of the origin $s=0$ and oriented
counterclockwise) and $\alpha$ is the angle (oriented clockwise)
between $u$ and the inner normal to $\bo$ in $q$. The Lebesgue \me\ on
$\ta \times S^1$ induces an $\ma$-invariant \me\ $\mu$ on $\ps$ which,
in the above coordinates, is described by $d\mu(s,\a) = c \, \cos\a \,
ds d\a$. The constant $c$ is customarily chosen so that $\mu$ is a
probability \me.

Let us indicate with $\si_0$ the set of all pairs $(s,\a) \in \ps$
where $s$ corresponds to a vertex of $\bo$ or $\a = \pm \pi/2$. The
set $\si_1 = \si_1^+$ introduced earlier is morally given by
``$\si_1^+ := \ma^{-1} \si_0$''. For historical reasons, this is
usually called the \emph{singularity set} of $\ma$, even though the
differential of $\ma$ is singular only at line elements resulting in
tangential hits. Analogously, for $n>1$, $\si_n^+ := \si_1^+ \cup
\ma^{-1} \si_1^+ \cup \cdots \cup \ma^{-n+1} \si_1^+$ is the set where
$\ma^n$ is not defined, which is called the singularity set of
$\ma^n$. We also introduce ``$\si_1^- := \ma \si_0$'' and, for $n>1$,
$\si_n^- := \si_1^- \cup \ma \si_1^- \cup \cdots \cup \ma^{n-1}
\si_1^-$. These are the singularity sets for the powers of the inverse
map $\ma^{-1}$. Lastly, $\si_\infty^+ := \bigcup_{n=1}^{\infty}
\si_n^+$, $\si_\infty^- := \bigcup_{n=1}^{\infty} \si_n^-$, and $\si
:= \si_\infty^+ \cup \si_\infty^-$.

Each $\si_n^\pm$ is the union of smooth curves whose endpoints lie
either on another such curve or on the \emph{generalized boundary} of
$\ps = [0,L] \times [-\pi/2,\pi/2]$, which is defined as the boundary
of $\ps$ plus all the vertical segments $s=s_i$, where $s_i$ is the
boundary coordinate of a vertex of $\bo$. If $L < \infty$, the number
of vertices is finite, and the curvature of $\ta$ is bounded, then
$\si_n^\pm$ comprises only a finite number of smooth curves.

Under the above assumptions, $\ma$ is a piecewise differentiable map
with singularities, of the type studied by Katok and Strelcyn in
\cite{ks}. Among their results is a suitable version of the Oseledec
Theorem which guarantees, for a.e.\ $(s,\a) =: x \in \ps$: 
\begin{enumerate}
\item A decomposition of the tangent space $T_x \ps$ into $E_x^+
  \oplus E_x^-$. These one-di\-men\-sion\-al spaces are
  dynamics-invariant in the sense that $(D \ma)_x E_x^\pm = E_{\ma
  x}^\pm$, where $(D \ma)_x$ denotes the differential of $\ma$ at
  $x$.

  \item The existence of the Lyapunov exponents $\lambda_\pm(x)$,
  defined as
  \begin{equation}
    \lambda_\pm(x) := \lim_{n \to +\infty} \frac1n \log \| (D \ma^n)_x
    v_\pm \|,
  \end{equation}
  with $v_\pm \in E_x^\pm$. Since $\mu$ is absolutely continuous
  w.r.t.\ the Lebesgue \me\ on $\ps$, then $\lambda_+ (x) = -\lambda_-
  (x)$. We adopt the convention that $\lambda_+ (x) \ge 0$.
\end{enumerate}
The \dsy\ is \hyp, by definition, if $\lambda_+ (x) > 0$ almost
everywhere. If the \sy\ is \erg\ too, then $\lambda_+ (x) =$ constant
$=: \lambda_+$.

\section{Geometrical optics and cone bundles}
\label{sec-cones}

In this section, which liberally draws from \cite{w2}, we recall the
basic tenets of the invariant cone technique for the \hyp ity
of planar \bi s (cf.\ also \cite{lw}), and prove a couple of results
that are specifically designed for our \sy s.

Given $x\in \ps$ and two linearly independent vectors $v_1, v_2 \in
T_x \ps$, we define the \emph{cone with boundaries $v_1, v_2$} as the
set 
\begin{equation}
  \label{cone}
  C(x) := \rset{av_1 + bv_2} {a,b \in \R,\ ab \ge 0 }. 
\end{equation}
If $C(x)$ is defined at every, or almost every, $x \in \ps$ and the
dependence on $x$ is measurable, we speak of $C \subset T\ps$ as a
measurable cone bundle.

A measurable cone bundle $C$ is said to be:
\begin{itemize}
\item \emph{invariant}, if $(D\ma)_x C(x) \subseteq C(\ma x)$ for
  $\mu$-a.e.\ $x$;

\item \emph{strictly invariant}, if $(D\ma)_x C(x) \subset C(\ma
  x)$ for $\mu$-a.e.\ $x$;

\item \emph{eventually strictly invariant}, if it is invariant
  and, for $\mu$-a.e.\ $x$, there exists $n(x) \in \Z^+$ such that $(D
  \ma^{n(x)})_x C(x) \subset C(\ma^{n(x)} x)$.
\end{itemize}

The next theorem was proved in \cite{w1}.

\begin{theorem}
  \label{thm-conhyp}
  Given a \bi\ map $\ma$ as described above, if there exists an
  eventually strictly invariant measurable cone bundle, then the
  Lyapunov exponent $\lambda_+(x)$ is positive for $\mu$-a.e.\ $x \in
  \ps$.
\end{theorem}

In \cite{w2} Wojtkowski reduces the invariance of a cone bundle to a
problem of geometrical optics concerning the behavior of a family (a
\emph{beam}) of nearby \tr ies. We present the main ideas.

To a tangent vector $v \in T_x \ps$ in phase space is naturally
associated a differentiable curve $\varphi: (-\eps, \eps)
\longrightarrow \ps$ such that $\varphi(0) = x$ and $\varphi'(0) = v$.
By construction, $\sigma \mapsto \varphi(\sigma)$ is uniquely
determined in linear approximation around $0$. Using the
representation of $\ps$ as a subset of $\ta \times S^1$, and the
notation $\varphi(\sigma) = (q(\sigma), u(\sigma)) \in \ta \times
S^1$, we construct the family of lines, or \emph{rays}, $l^+ (\sigma)
:= \lset{q(\sigma) + r u(\sigma)} {r \in \R}$. Also, denoting by
$u^-(\sigma)$ the outward-pointing, precollisional vector of
$u(\sigma)$ at $q(\sigma) \in \bo$, we define $l^- (\sigma) :=
\lset{q(\sigma) + r u^- (\sigma)} {r \in \R}$.

In first approximation, that is, when $\eps \to 0^+$, the now
infinitesimal beam of rays \emph{focuses} in a point, which means that
all rays, up to adjustments of order $\eps$ in $(q(\sigma),
u(\sigma))$, have a common intersection. We consider the case too
where the common intersection is at infinity. This \emph{focal point}
is clearly a \fn\ of $v$ only: it is denoted $F^+(v)$ for the family
$\{ l^+ (\sigma) \}$ and $F^-(v)$ for the family $\{ l^- (\sigma) \}$.
Let us call $f^\pm (v)$ the signed distances, along $l^\pm(0)$,
between $F^\pm(v)$ and $q_0 = q(0)$ ($l^\pm(\sigma)$ has the
orientation induced by the parameter $r \in \R$, that is, outward for
$l^- (\sigma)$ and inward for $l^+ (\sigma)$, relative to $\ta$). In
the remainder, we will omit the dependence of $v$ from all the
notation whenever there is no ambiguity.  Indicated by $(ds, d\a)$ the
components of $0 \ne v \in T \ps_{(s_0, \a_0)}$ in the natural basis
$\{ \partial / \partial s,\partial / \partial\a \}$, one has
\begin{equation}
  f^{\pm} = \left\{
  \begin{array}{lll}
    \ds \frac{\cos \a_0} {\pm k(s_0) - \frac{d\a}{ds}}, && 
      \mbox{if } ds \ne 0; \vspace{6pt} \\
    0, && \mbox{if } ds = 0.
  \end{array} 
  \right.
  \label{fpm}
\end{equation}
Here $k(s)$ denotes the curvature of $\bo$ at the point of coordinate
$s$ (as specified in the introduction, the curvature is taken positive
at focusing points of the boundary, and negative at dispersing
points). The formula (\ref{fpm}) is derived, e.g., in \cite{w2}.

It is easy to see that $f^\pm$ are projective coordinates of $T_x
\ps$. Hence any cone of the type (\ref{cone}) can be described by a
closed interval in the coordinate $f^+ \in \overline{\R}$, where
$\overline{\R} := \R \cup \{ \infty \}$ is the compactification of
$\R$.  Henceforth, for simplicity, we will drop the subscripts from
the coordinates $(s_0, \a_0)$ of the collision pair. Also, we will use
the imprecise terminology `the point $s \in \bo$' to mean `the point
in $\bo$ of coordinate $s$'.  The next lemma is known in optics as the
mirror equation \cite{w2, cm}.

\begin{lemma}
  For an infinitesimal beam of \tr ies colliding around the point $s
  \in \bo$ with reflection angles around $\a$,
  \begin{displaymath}
    -\frac1{f^{-}} + \frac1{f^{+}} = \frac{2k(s)} {\cos \a}.
  \end{displaymath}
\label{lem-mirror-eq}
\end{lemma}

We now present a visual description of the cone $C(x) = C(s,\a)$ on
the configuration plane containing $\ta$.  For $s \in \bo$ and
$\beta>0$, denote by $D_\beta (s)$ the closed disc of radius $1 /
|\beta k(s)|$ tangent to $\bo$ in $s$ on the internal side of $\ta$.
Analogously, for $\beta<0$, let $D_\beta (s)$ be the closed disc of
radius $1 / |\beta k(s)|$ tangent to $\bo$ in $s$ on the external side
of $\ta$.  Consider also the two closed halfplanes delimited by
$t(s)$, the tangent line to $\bo$ in $s$: let $D_{0+}(s)$ denote the
internal halfplane, relative to $\ta$, and $D_{0-}(s)$ the external
one. See Fig.~\ref{fig-dbetas}. The interior of $D_\beta (s)$ is
indicated with $D_\beta^\circ (s)$.

\fig{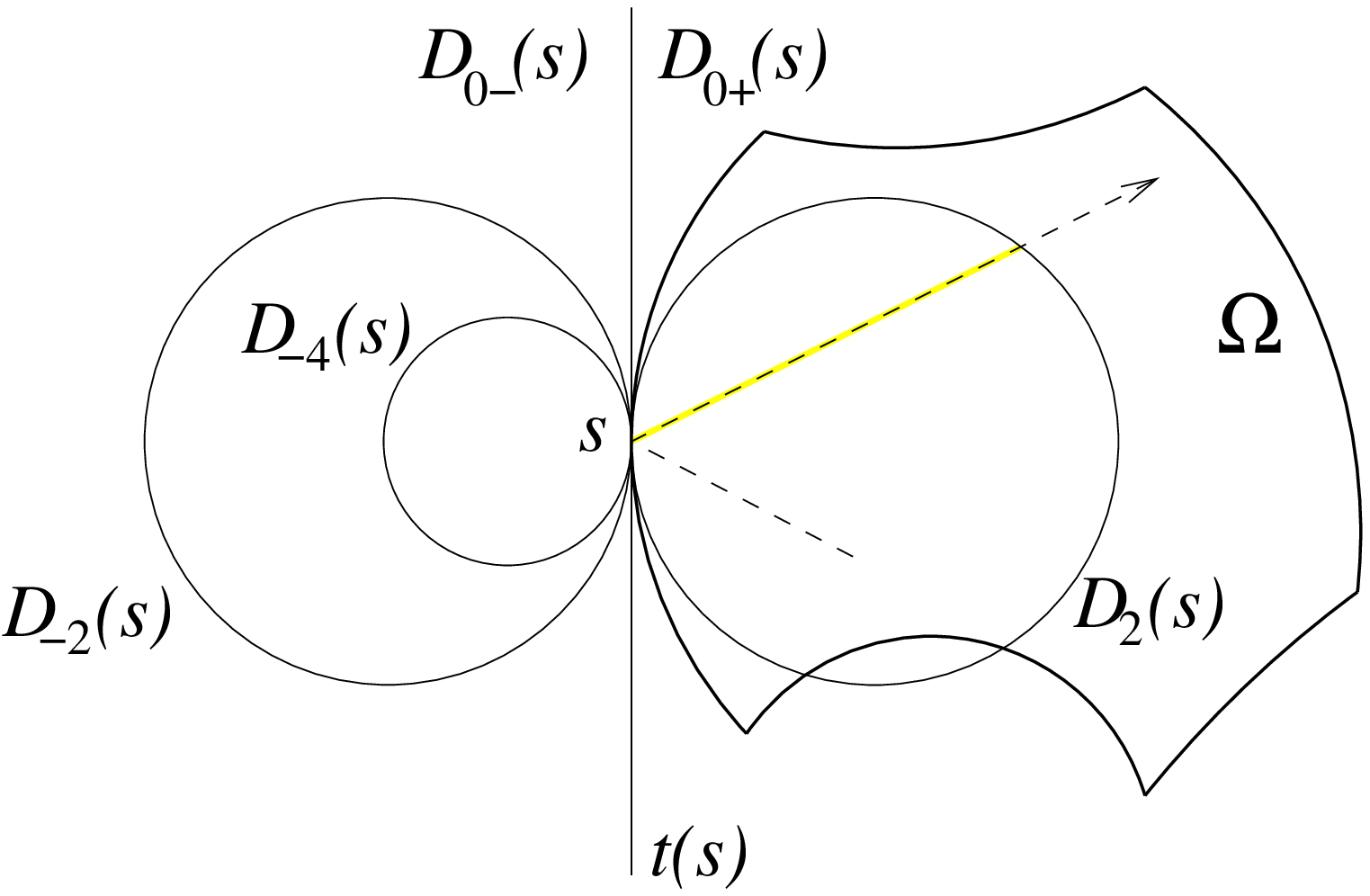} {8cm} {The tangent line $t(s)$ and some discs
  $D_\beta(s)$. The yellow part of the \tr y is the locus of the focal
  points $F^+$ corresponding to a certain cone.}

\begin{lemma} 
  Given a cone $C(s,\a)$ of the type \emph{(\ref{cone})}, $v \in
  C(s,\a)$ corresponds to $F^+(v) \in l^+(0) \cap D$, where $D
  \subset \R^2$ is one of the following sets:
  \begin{itemize}
  \item[(a)] $D = D_{\beta_1} (s)$;

  \item[(b)] $D = D_{\beta_1} (s) \setminus D_{\beta_2}^\circ (s)$,
    with $|\beta_1| < |\beta_2|$;

  \item[(c)] $D = D_{\beta_1} (s) \cup D_{\beta_2} (s)$, with $\beta_1 
    \ge 0$ and $\beta_2 \le 0$;

  \item[(d)] $D = \R^2 \setminus (D_{\beta_1}^\circ (s) \cup
    D_{\beta_2}^\circ (s) \cup \{s\} )$, with $\beta_1 \ge 0$ and
    $\beta_2 \le 0$.
  \end{itemize}
  Moreover,
  \begin{displaymath}
    F^+(v) \in \partial D_\beta (s) \setminus \{ s \} \quad 
    \Longleftrightarrow \quad f^+(v) = \frac{2 \cos \a} {\beta 
    |k(s)|}.
  \end{displaymath}
  \label{lem-fplus}
\end{lemma} 

\proof By construction $F^+ = F^+(v) \in l^+(0)$.  Since $f^+$ is a
coordinate on $l^+(0)$, a closed interval in the projectivized $f^+
\in \overline{\R}$ corresponds, on $l^+(0)$, to either a closed
segment or a closed halfline or the union of two disjoint closed
halflines. Cases \emph{(a)}-\emph{(d)} cover all possibilities.

The second statement, for $\beta>0$, comes from elementary
trigonometry (see Fig.~\ref{fig-dbetas}), and it trivially extends to
the case $\beta<0$ as well.  
\qed

The reason why, in Lemma \ref{lem-fplus}, we chose such peculiar sets
$D$ to cut a (projective) closed segment on $l^+(0)$, upon
intersection, will be made clear by the next lemma. In particular, we
will see that describing the cones in terms of the discs $D_\beta(s)$
will eliminate the dependence on $\a$ in the mirror equation of Lemma
\ref{lem-mirror-eq}.

\begin{lemma} 
  For infinitesimal beam of trajectories colliding around $s \in \bo$,
  $F^{-} \in \partial D_\beta (s)$ \iff $F^{+} \in \partial D_{\beta'}
  (s)$, where
  \begin{displaymath}
    \beta' = 4 \, \mathrm{sgn}(k(s)) - \beta
  \end{displaymath}
  (with the understanding that $F^\pm \in \partial D_{0\pm}$ means
  $F^\pm \in \{ s, \infty \}$).
  \label{lem-betap}
\end{lemma}

\proof Let $\a$ be the angle of reflection (and thus of incidence) of
the \tr y we are perturbing. Disregarding the case $F^+ = F^- = s$, we
know from Lemma \ref{lem-fplus} that $F^+ \in \partial D_{\beta'} (s)$
corresponds to $f^+ = 2 \cos\a / (\beta' |k(s)|)$. Also, $F^- \in
\partial D_\beta (s)$ is equivalent to $f^- = -2 \cos\a / (\beta
|k(s)|)$ (the minus sign is needed because a focal point $F^-$ lying
on the internal halfplane $D_{0+} (s)$ corresponds to a negative $f^-$
along $l^-(0)$, and viceversa). Direct substitution into Lemma
\ref{lem-mirror-eq} yields
\begin{equation}
  \frac{\beta |k(s)|} {2\cos \a} + \frac{\beta' |k(s)|} {2 \cos\a} =  
  \frac{2k(s)} {\cos\a},
\end{equation}
whence the assertion.
\qed

With the tools of Section \ref{sec-cones}, the problem of the cone
invariance along a given \tr y can be reduced to the study of the
focal points of one-parameter perturbations of that \tr y.

We single out the information that we need for our forthcoming proofs.

\begin{proposition}
  For an infinitesimal beam of \tr ies colliding around $s$ we have
  the following:
  If $s$ belongs to a focusing component of $\bo$, i.e., $k(s)>0$, then:
  \begin{eqnarray*}
    F^\mp \in D_4(s) & \Longleftrightarrow & F^\pm \in D_{0-}(s); \\
    F^\mp \in D_2(s) \setminus D_4^\circ (s) & \Longleftrightarrow & 
    F^\pm \in D_{0+}(s) \setminus D_2^\circ (s).
  \end{eqnarray*}
  If $s$ belongs to a dispersing component of $\bo$, i.e., $k(s)<0$,
  then
  \begin{eqnarray*}
    F^\mp \in D_{-4}(s) & \Longleftrightarrow & F^\pm \in D_{0+}(s); \\
    F^\mp \in D_{-2}(s) \setminus D_{-4}^\circ (s) & \Longleftrightarrow
    & F^\pm \in D_{0-}(s) \setminus D_{-2}^\circ (s).
  \end{eqnarray*}
  Analogous equivalences hold for the interior of such cones. The
  situation is illustrated in Fig.~\ref{fig-p-dbeta}.
  \label{prop-d-beta}
\end{proposition}

\proof We only prove the first statement, the other ones
being completely analogous. Once again, we disregard the easy case
$F^+ = F^- = s$.  We have $F^- \in D_4(s)$ $\Leftrightarrow$ $F^- \in
\partial D_\beta(s)$, for $\beta \in [4, +\infty)$ $\Leftrightarrow$
(by Lemma \ref{lem-betap}) $F^+ \in \partial D_{\beta'} (s)$, for
$\beta' \in (-\infty, 0]$ $\Leftrightarrow$ $F^- \in D_{0-}
(s)$. Clearly, nothing changes if we swap $F^-$ and $F^+$.
\qed

\fig{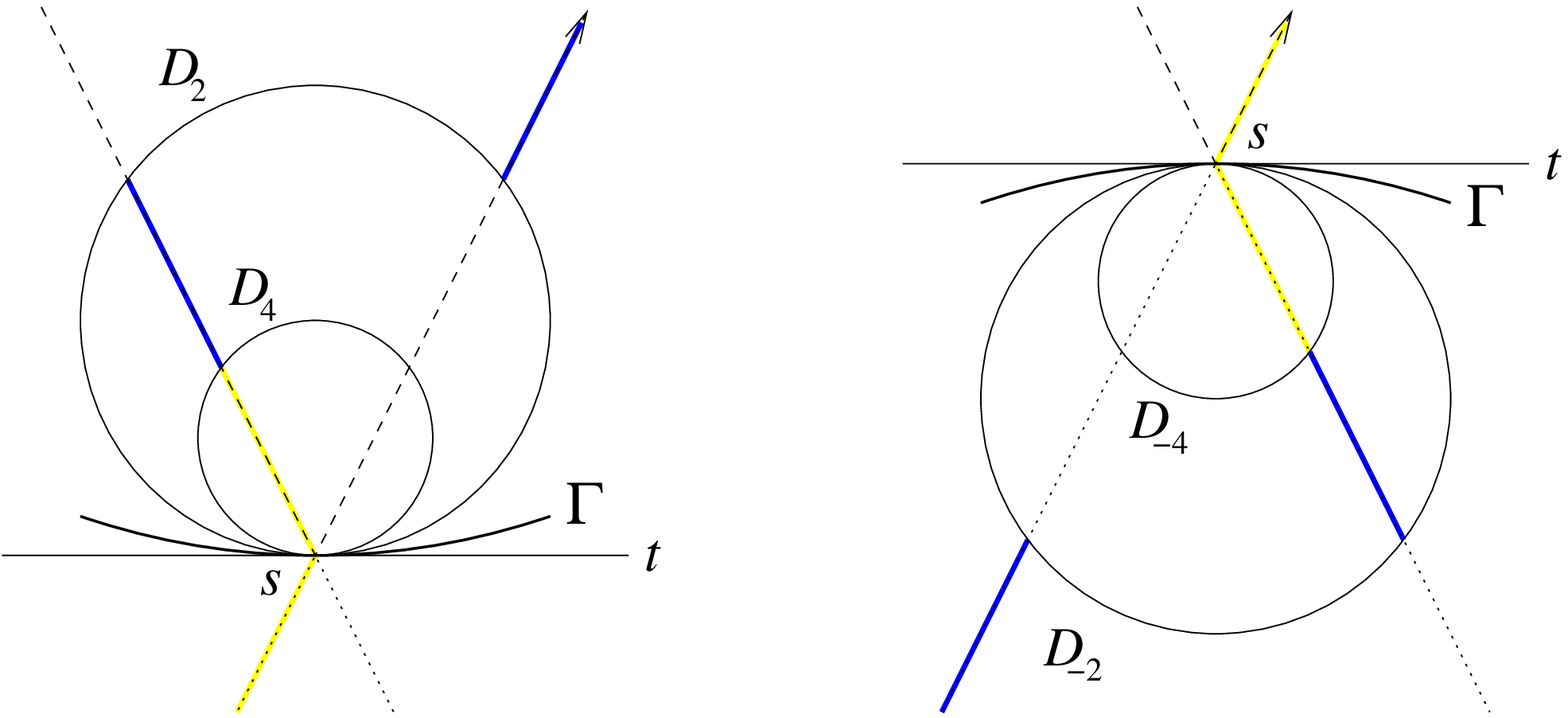} {13.5cm} {A geometric representation of Proposition
  \ref{prop-d-beta}. The left picture represents the first two
  statements (focusing border); the right picture represents the last
  two statements (dispersing border). Yellow/blue sets of focal points
  $F^-$ are mapped into yellow/blue sets of focal points $F^+$. The
  dependence on $s$ in the notation has been omitted.}

\section{Hyperbolicity}
\label{sec-hyp}

Fig.~\ref{fig-t1} shows the \bi\ table we are mainly interested in for
the rest of the paper. We refer to it for the definition of the
quantities $l, h > 0$. The three dispersing components of the boundary
$\bo$ are circular arcs of curvature $k_d \in (-\sqrt{2},0)$.  Their
union is denoted $\bo_d$. The focusing component is a circular arc of
curvature $k_f > 0$ and is denoted $\bo_f$. The remining, flat, part
of the boundary is denoted $\bo_s$. The two rectangular portions of
$\ta$ which $\bo_s$ almost delimits will be referred to as \emph{the
strips}, or \emph{the corridors}, or whatever one's fancy suggests
each time.

\fig{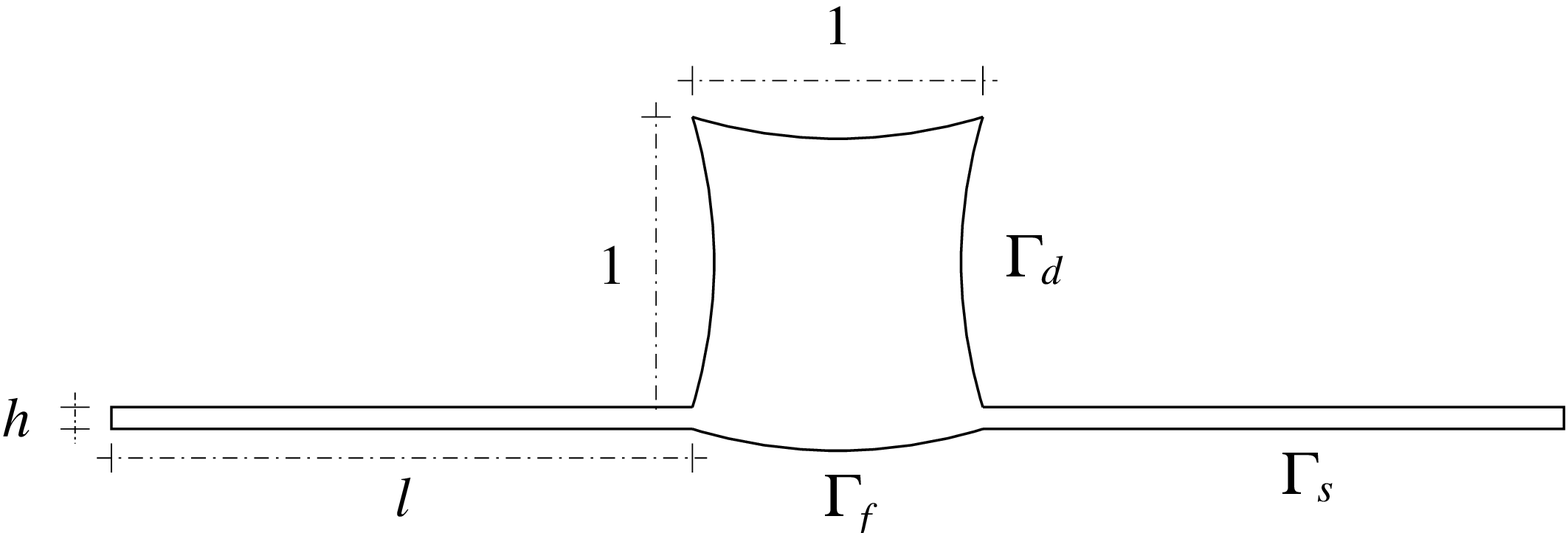} {12.4cm} {The definition of the table $\ta$.}

The geometric constants $l, h, k_f, k_d$ are chosen via the following
procedure. Keep in mind that we are interested in small values of
$k_f$ (see the Introduction) and $h$ (see Section \ref{sec-conf}).
One starts by fixing arbitrary values of $k_d$ and $h$.  Then $k_f$ is
determined by a geometric condition that we presently describe, with
the help of Fig.~\ref{fig-iss}. For $s' \in \bo_d$ and $s'' \in
\bo_f$, consider the straight line passing through $s'$ and $s''$, and
let $I(s', s'')$ be its intersection with the disc $D_{-2}(s')$. The
curvature $k_f$ must be so small that
\begin{equation}
  \label{C1}
  \forall s' \in \bo_d, \ \forall s'' \in \bo_f, \quad I(s', s'') 
  \subset D_4(s''). 
\end{equation}
Finally, $l$ is chosen such that
\begin{equation}
  \label{C2}
  l \ge \frac1 {k_f}
\end{equation}

\fig{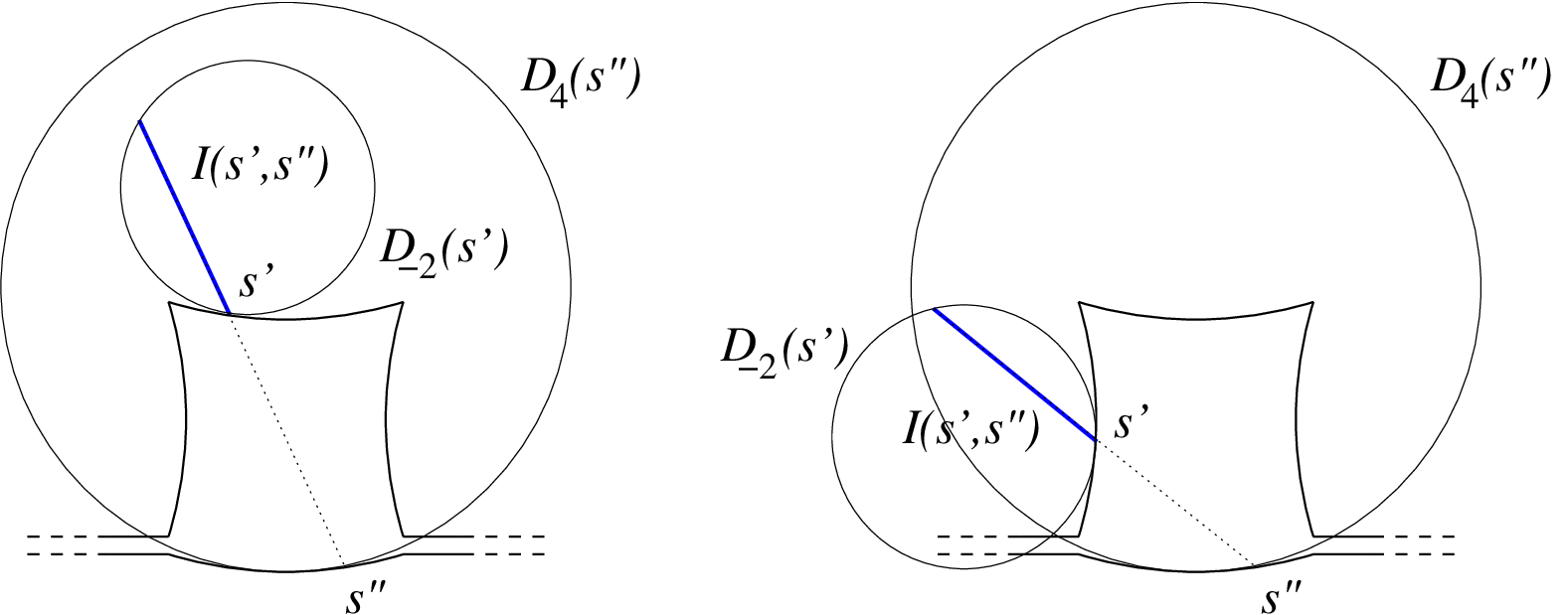} {13cm} {Condition (\ref{C1}) for two different choices
  of $s'$.}

\begin{remark}
  Condition (\ref{C1}) excludes sufficient separation between the
  boundary components as per the standard theory of Wojtkowski,
  Markarian, Donnay and Bunimovich, which is summed up, e.g., in
  \cite[Thm.~9.19]{cm}. The hypotheses of that theorem are evidently
  violated as (\ref{C1}) implies in particular that $D_4(s'')$
  contains large portions of $\bo_d$, for all $s \in \bo_f$.
  \label{rk-sep}
\end{remark}

We are now going to prove the \hyp ity of this \bi\ \sy\ via Theorem
\ref{thm-conhyp}. However, we will not use exactly the Poincar\'e
section that we have introduced in Sections \ref{sec-prel} and
\ref{sec-cones}, but a similar section that neglects the hits on the
flat boundary component $\bo_s$. This is standard procedure in the
theory of \hyp\ \bi s as it is basic fact that the collisions against
a flat boundary do not change the \hyp\ features of a beam of \tr ies.
(One easy way to see this is to \emph{unfold} the \bi\ along a given
\tr y: every time the material point hits a flat side we pretend that
it continues its precollisional rectilinear motion, but we reflect the
table around that flat side; apart from this rigid motion of the \bi\
table, nothing changes for the \tr y or any of its infinitesimal
perturbations.)

Let us denote $\bar{\bo} := \bo_f \cup \bo_d$. With the usual abuse of
notation, whereby a point $q \in \bo$ is identified with its arclength
coordinate $s$, we define $\ps := \bar{\bo} \times [-\pi/2, \pi/2]$,
whose elements we call $(s,\a)$ or $x$. Clearly $\ps$ is a global
cross section for the flow. Let $\ma: \ps \longrightarrow \ps$ be its
first-return map.

For any $x = (s, \a) \in \ps$ and $n \in \Z$, denote $x_n := (s_n,
\a_n) := \ma^n x$ and let $\tau_n$ be the length of the portion of the
\tr y (equivalently, the time) between the collisions at $s_n$ and
$s_{n+1}$ (notice that there can be an arbitrary number of collisions
against $\bo_s$ between $s_n$ and $s_{n+1}$). Also, let $k_n :=
k(s_n)$ indicate the curvature of $\bo$ in $s_n$. Analogously, given
$v \in T_x \ps$, denote $v_n := (D\ma^n)_x v$. The infinitesimal
beam of \tr ies determined by $v_n$ (and thus by $v$) around $(s_n,
\a_n)$ will have pre- and postcollisional foci denoted, respectively,
$F_n^- := F^-(v_n)$ and $F_n^+ := F^+(v_n)$. The corresponding signed
distances along the pre- and postcollisional lines are indicated with
$f_n^-$ and $f_n^+$. The following facts are obvious:
\begin{eqnarray}
  && F_n^- = F_{n-1}^+, \\
  && f_n^- = -( \tau_{n-1} - f_{n-1}^+).
\end{eqnarray}
For the sake of the notation, let us drop all subscripts 0 and
write $k := k_0$, $F^+ := F_0^+$, and so on.  

For any $x \in \ps$, we introduce the following three cones in $T_x
\ps$:

\begin{itemize}
\item $C_0 (x)$ is the set of all tangent vectors whose correspondent
  family of rays focuses in linear approximation inside $D_{-2}(s)$.
  Using the focal distance $f^+$,
  \begin{equation}
    C_0 (x) := \rset{v \in T_x \ps} {-\frac{\cos \alpha} {|k|} \le
    f^+(v) \le 0}.
  \end{equation}

\item $C_1 (x)$ is the set of all tangent vectors whose correspondent
  family of rays focuses in linear approximation inside $D_{0-}(s)$,
  i.e., all the divergent families of rays. In projective terms,
  \begin{equation}
    C_1 (x) := \rset{v \in T_x \ps} {-\infty <  f^+(v) \le 0}.
  \end{equation}

\item $C_2 (x)$ is the set of all tangent vectors whose correspondent
  family of rays focuses in linear approximation inside $D_2 (s)
  \setminus D_4^\circ (s)$, i.e.,
  \begin{equation}
    C_2 (x) := \rset{v \in T_x \ps} {\frac{\cos \alpha} {2|k|} \le 
    f^+(v) \le \frac{\cos \alpha} {|k|}}.
  \end{equation}
\end{itemize}

We use the above cones to define piecewise an invariant cone bundle $C
:= \{ C(x) \}_x$. For each $x = (s, \a)$, the choice $C(x) := C_i (x)$
will depend on $s$, $s_{-1}$, and what happens to the \tr y between
the collisions at $s_{-1}$ and $s$.

\begin{itemize}
\item[(A)] If \underline{$s \in \bo_d$}, set $C(x) := C_0 (x)$.

\item[(B)] If \underline{$s \in \bo_f$}, there are two subcases:
  \begin{itemize}
  \item[(B.1)] If \underline{$s_{-1} \in \bo_f$}, set $C(x) := C_2 (x)$.

  \item[(B.2)] If \underline{$s_{-1} \in \bo_d$}, there are two further 
    subcases, depending on whether the piece of \tr y between $s_{-1}$ 
    and $s$ has collisions with $\bo_s$:
    \begin{itemize}
    \item[(B.2.1)] \underline{No collisions with $\bo_s$} between 
      $s_{-1}$ and $s$: Set $C(x) := C_1 (x)$.

    \item[(B.2.2)] \underline{At least one collision with $\bo_s$} 
      between $s_{-1}$ and $s$: Set $C(x) := C_2 (x)$.
    \end{itemize}
  \end{itemize}
\end{itemize}

Clearly $C(x)$ is a measurable \fn\ of $x$.

\begin{theorem}
  \label{thm-hyp}
  The cone bundle $C$ just defined is eventually strictly invariant
  relative to the map $\ma$.
\end{theorem}

\proof We check that $v \in C(x)$ implies $v_1 \in C(x_1)$ for all the
possible cases $C(x) = C_i(x)$, $C(x_1) = C_j(x_1)$ $(i,j \in \{ 0,1,2
\})$.

\begin{itemize}

\item[(I)] \underline{$s, s_1 \in \bo_d$}.\ In this case $C(x) =
  C_0(x)$, $C(x_1) = C_0(x_1)$. $v \in C_0(x)$ implies $F^+ \in
  D_{-2}(s)$, hence $F_1^- = F^+ \in D_{0+}^\circ (s_1)$. By
  Proposition \ref{prop-d-beta}, $F_1^+ \in D_{-4}^\circ (s_1) \subset
  D_{-2}^\circ (s_1)$.  This is equivalent to $v_1 \in C_0^\circ
  (x_1)$---where $C^\circ(x)$ represents the interior of $C(x)$ in
  $T_x \ps$. We have thus proved strict invariance for this type of
  collision.

\item[(II)] \underline{$s \in \bo_d$, $s_1 \in \bo_f$}.\ Here $C(x) =
  C_0(x)$ but the cone $C(x_1)$ may take two different forms. We 
  separately check both cases.

  \begin{itemize}
  \item[(II.1)] There are no collisions with $\bo_s$ between $s$ and
    $s_1$. Then $C(x_1) = C_1(x_1)$. For $v \in C_0(x)$ we have, by
    condition (\ref{C1}), $F_1^- = F^+ \in D_4 (s_1)$. Proposition
    \ref{prop-d-beta} implies that $F_1^+ \in D_{0-} (s_1)$, that is,
    $v_1 \in C_1 (x_1)$. In this case the invariance is not
    necessarily strict.

  \item[(II.2)] There are collisions with $\bo_s$ between $s$ and
    $s_1$, that is, the material point enters a strip before colliding
    at $s_1$. In this case $C(x_1) = C_2(x_1)$. Since the material
    point has to travel all the way to the end of the strip and bounce
    back, $\tau > 2l > 2/k_f$, having used condition (\ref{C2}). For
    $v \in C_0(x)$, $f^+ \le 0$, hence $f_1^- = -\tau + f^+ < -1/k_f$.
    Equivalently, $F_1^- \in D_{0+}(s_1) \setminus D_2(s_1)$. By
    Proposition \ref{prop-d-beta}, $F_1^+ \in D_2^\circ (s_1)
    \setminus D_4 (s_1)$, i.e., $v_1 \in C_2^\circ (x_1)$.
  \end{itemize}

\item[(III)] \underline{$s \in \bo_f$, $s_1 \in \bo_d$}.\ Here $C(x_1)
  = C_0(x_1)$ and we have two subcases on $C(x)$.

  \begin{itemize}
  \item[(III.1)] $C(x) = C_1(x)$. In this case $v \in C(x)$ is
    equivalent to $f^+ \le 0$. Hence $f_1^- < 0$ and $F_1^- \in
    D_{0+}^\circ (s_1)$. Therefore (Proposition \ref{prop-d-beta})
    $F_1^+ \in D_{-4}^\circ (s_1) \subset D_{-2}^\circ (s_1)$. Namely
    $v_1 \in C_0^\circ (x_1)$.

  \item[(III.2)] $C(x) = C_2(x)$. So $v \in C(x)$ means that $F^+ =
    F_1^- \in D_2 (s) \setminus D_4^\circ (s)$. We consider two
    possible types of \tr ies:

    \begin{itemize}
    \item[(III.2.1)] There are no collisions with $\bo_s$ between $s$
      and $s_1$. By (\ref{C1}), $F_1^- \in D_{0-}^\circ (s_1)
      \setminus D_{-2}^\circ (s_1)$. Hence $F_1^+ \in D_{-2} (s_1)$.

    \item[(III.2.2)] There are collisions with $\bo_s$ between $s$ and
      $s_1$. As in case (II.2), $\tau > 2/k_f$ and $f^+ \le (\cos\a) /
      k_f < 0$. Thus, $f_1^- < 0$, that is, $F_1^- \in D_{0+}^\circ
      (s_1)$. Finally, $F_1^+ \in D_{-4}^\circ (s_1) \subset
      D_{-2}^\circ (s_1)$.
    \end{itemize}
  \end{itemize}

\item[(IV)] \underline{$s, s_1 \in \bo_f$}.\ Definition (B.1) ensures
  that $C(x_1) = C_2(x_1)$. Let us branch out in two subcases
  depending on $C(x)$.

  \begin{itemize}
  \item[(IV.1)] $C(x) = C_1(x)$. As in case (III.1), $v \in C(x)$
    implies that $f^+ \le 0$. Since, by construction of our
    cross section, there can be no collisions with $\bo_d$ in the
    piece of \tr y between $s$ and $s_1$, there are only two
    possibilities: either the particle enters and exits a strip, and
    thus $\tau > 2/k_f$; or that piece of \tr y is a chord of the arc
    $\bo_f$, and thus $\tau = 2 (\cos\a) /k_f$. In either case, $\tau
    > (\cos\a) /k_f$ and $f_1^- < -(\cos\a) /k_f$, which means that
    $F_1^+ \in D_{0+}^\circ (s_1) \setminus D_2 (s_1)$. By Proposition
    \ref{prop-d-beta}, $F_1^+ \in D_2^\circ (s_1) \setminus D_4
    (s_1)$, that is, $v_1 \in C_2^\circ (x_1)$.

  \item[(IV.2)] $C(x) = C_2(x)$. The hypothesis $v \in C(x)$ reads
    $(\cos\a) / 2k_f \le f^+ \le (\cos\a) / k_f$. Once again, there are
    two further subcases:

    \begin{itemize}
    \item[(IV.2.1)] There are no collisions with $\bo_s$ between $s$
      and $s_1$. In this case, cf.\ (IV.1), the \tr y between $s$ and
      $s_1$ is a chord of $\bo_f$ and $\tau = 2 (\cos\a) /k_f$.
      Therefore $f_1^- = -\tau + f^+ \le -(\cos\a) /k_f$, which
      implies $F_1^- \in D_{0+} (s_1) \setminus D_2^\circ (s_1)$. This
      yields $F_1^+ \in D_2 (s_1) \setminus D_4^\circ (s_1)$, namely
      $v_1 \in C_2(x_1)$. 

    \item[(IV.2.2)] There are collisions with $\bo_s$ between $s$ and
      $s_1$. $f^+$ and $\tau$ are exactly as in case (III.2.2).
      Refining the estimate that is written there, $f_1^- < -1/k_f <
      -(\cos\a) / k_f$, that is, $F_1^- \in D_{0+}^\circ (s_1)
      \setminus D_2 (s_1)$. This gives $F_1^+ \in D_2^\circ (s_1)
      \setminus D_4 (s_1)$.
    \end{itemize}

  \end{itemize}

\end{itemize}

In order to show that $C$ is eventually strict invariant \emph{almost
everywhere}, we notice that there are only three cases above in
which the cone invariance is not strict, namely (II.1), (III.2.1), and
(IV.2.1).

In both (II.1) and (III.2.1), nonstrictness can only occur when the
external endpoint of $I(s',s'')$ lies on $D_4 (s'')$ and $s = s'$,
$s_1 = s''$, or viceversa---cf.\ (\ref{C1}) and Fig.~\ref{fig-iss}. It
is not hard to realize that this situation can only occur for finitely
many pairs $(s',s'')$ (at least when the table is optimized, see
(\ref{C3}) and Fig.~\ref{fig-ho}, there are only two such pairs).

As concerns (IV.2.1), we realize that there can only be a finite
number of consecutive collisions of that type, because each such piece
of \tr y is a chord of $\bo_f$ of constant length ($\tau = \tau_1$),
but $\bo_f$ is smaller than a semicircle.
\qed

\section{Confining the table to a bounded region}
\label{sec-conf}

In the previous section the table $\ta$ was constructed starting with
two values for $h$ and $k_d$, which determined an upper bound on the
choice of $k_f$, via (\ref{C1}), which in turn determined a lower
bound on the choice of $l$, via (\ref{C2}). The latter condition, in
particular, forced the area of $\ta$ to diverge, as smaller and smaller
values are chosen for $k_f$.

Now we want to optimize, that is, minimize, the area of the table and
to do so we change the order in which its geometric parameters are
chosen. Given $k_d < 0$ and $k_f$ sufficiently small, we define the
\emph{optimal height} and the \emph{optimal length} of the strips,
respectively, as:
\begin{eqnarray}
  \label{C3}
  && h_o := h_o (k_d, k_f) := \min \rset{h} {\forall s' \in \bo_d,  
  \forall s'' \in \bo_f, \  I(s', s'') \subset D_4(s'')} ; \\
  \label{C4}
  && l_o := l_o (k_f) = k_f^{-1} .
\end{eqnarray}
These definitions are well posed, in the sense that a table can be
constructed with $h = h_o$ and $l = l_o$. We call it the \emph{optimal
table} and we think of it as a \fn\ of $k_f$ ($k_d$ is considered
fixed once and for all). The optimal table is \hyp\ by Theorem
\ref{thm-hyp}. The next proposition shows that, as $k_f \to 0$, the
area of the optimal table is bounded above. (In what follows, the
notation $a \sim b$ means that $a = a(k_f)$, $b = b(k_f)$ and, as $k_f
\to 0$, $|a/b|$ is bounded away from $0$ and $\infty$.)

\begin{proposition}
  \label{prop-ho}
  As $k_f \to 0$, $h_o(k_f) \sim k_f$.
\end{proposition}

\proof Since $k_f \to 0$ and $k_d$ is fixed, we may assume that, given
any $s'' \in \bo_f$, $D_4(s'')$ easily contains $D_{-2}(s')$, for all
$s'$ in the upper component of $\bo_d$ (left picture in
Fig.~\ref{fig-iss}).  

For $s'$ belonging to the lateral components of $\bo_d$, it is not
hard to realize that the worst-case scenario is the one depicted in
Fig.~\ref{fig-ho} (or the specular situation w.r.t.\ the axis of
symmetry of $\ta$): First of all, if $s''$ moves to the left and/or
$s'$ moves upward, $I(s',s'')$ will move towards the interior of
$D_4(s'')$, so that (\ref{C1}) is always verified.  Secondly, setting
$h_o$ to be the $h$ displayed there, one clearly sees that for $h \ge
h_o$ (\ref{C1}) is verified, while for $h < h_o$ it is not.

\fig{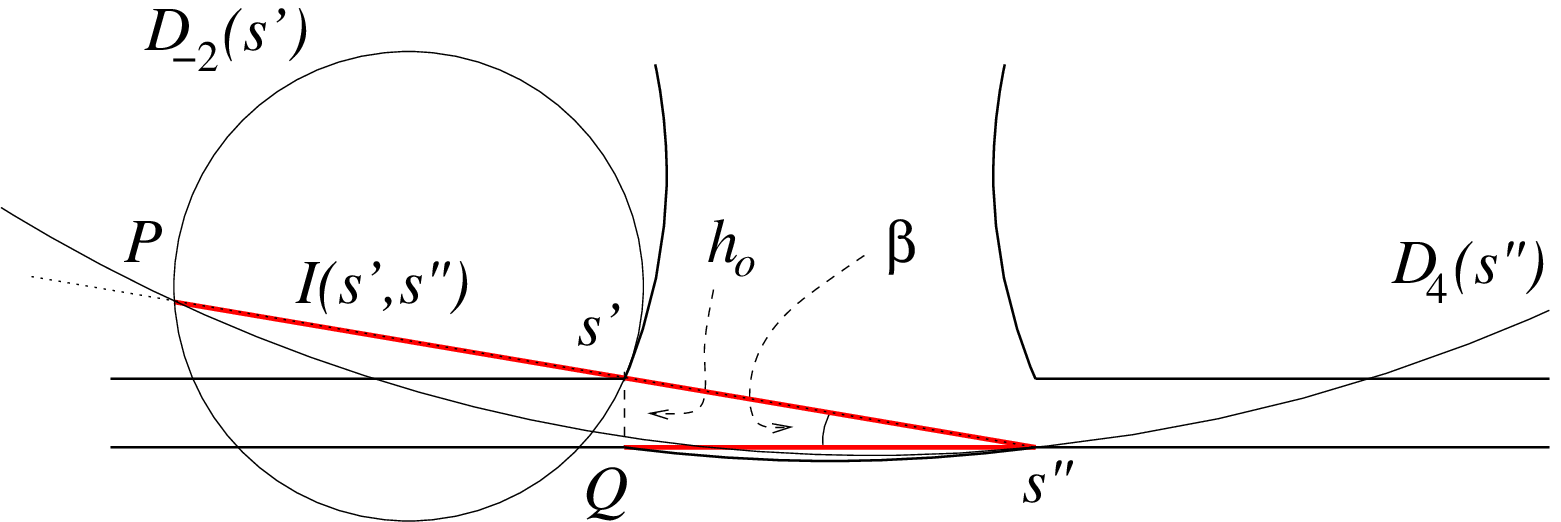} {10cm} {Finding $h_o$, cf.\ Proposition \ref{prop-ho}.}

Referring to the notation of Fig.~\ref{fig-ho}, we see that $h_o =
\tan \beta$ where $\beta$ is the angle between the two chords $s''P$
and $s''Q$ of $\partial D_4(s'')$. Recalling that, in a circle of
radius $r$, the relation between the length $\ell$ of a chord and the
angle $\theta$ it makes with the tangent to the circle at each of its
endpoints is $\ell = 2r \sin \theta$, we have
\begin{equation}
  \beta = \arcsin \left( \frac{k_f \, c}2 \right) - \arcsin \left(
  \frac{k_f}2 \right) \sim k_f, \quad \mbox{as } k_f \to 0.
\end{equation}
In the above $c$ is the length of $s''P$, for which it holds $1 < c <
2 + 2k_d^{-1}$. This ends the proof since $h_o \sim \beta$.
\qed

From a technical point of view, Proposition \ref{prop-ho} is a
consequence of the fact that $\bo_f$ fails to act as a perturbation of
a semidispersing component only for a few \tr ies, whose corresponding
beams need to be defocused by visiting the long strips. As $k_f \to
0$, this phenomenon concerns fewer and fewer \tr ies, but its fix
requires more and more space.  Proposition \ref{prop-ho} tells us that
the trade-off between the two effects balances out.

If a \hyp\ \bi\ table with a flatter and flatter focusing component
need not become bigger and bigger in terms of area, one might hope
that it need not in terms of diameter, either. In our particular
table, one would like to redesign the strips so that their area is
better placed in the plane and can be included in a fixed compact
region. In the remainder of the section we show that this is possible,
for example by bending the strips around the bulk of the \bi\
(see Fig.~\ref{fig-t2intro}).

\fig{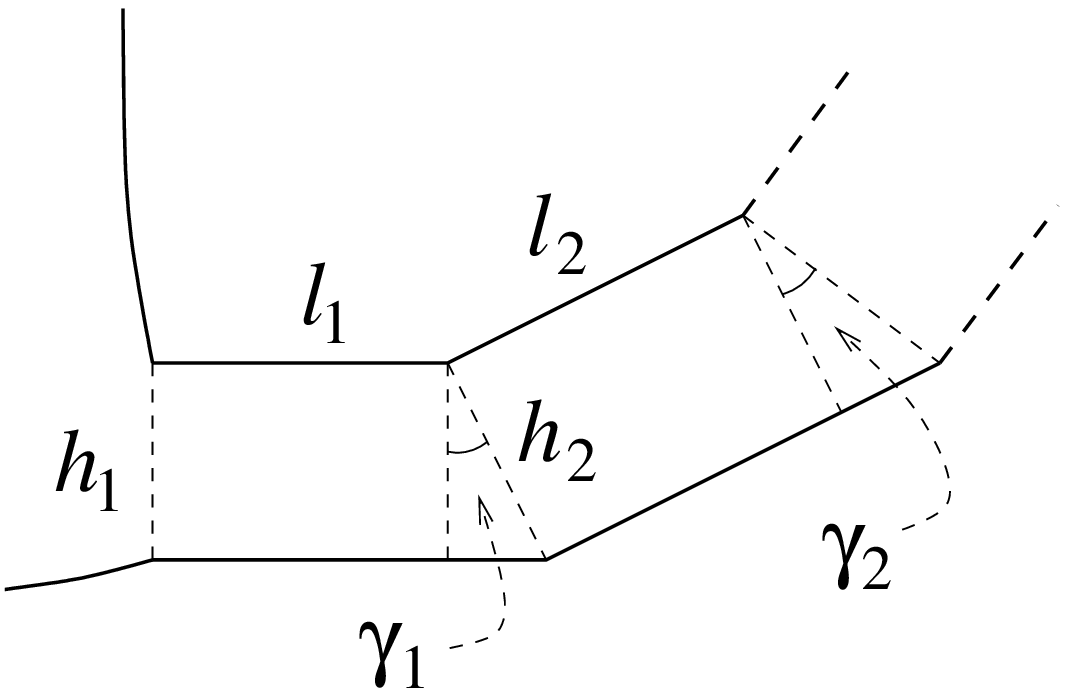} {5cm} {Construction of a spiral as a union of
  trapezoids.}

Let us describe this construction with the help of
Fig.~\ref{fig-folding}. Substitute each strip of $\ta$ with a
polygonal modification given by the union of $N$ adjacent right
trapezoids $T_1, T_2, \ldots, T_N$, where $N$ will be specified later
depending on $k_f$. $T_1$ is placed so that its shorter leg coincides
with the opening towards the bulk of $\ta$: its height is then $h_1 :=
h \ge h_o$. The length of the shorter base is denoted $l_1$ and the
two nonright angles are denoted $\pi/2 + \gamma_1$ and $\pi/2 -
\gamma_1$, with $0 < \gamma_1 < \pi/2$. This causes the longer leg to
measure $h_2 := h_1 / \cos \gamma_1$. The longer leg of $T_1$ is then
used as the shorter leg of the next trapezoid, $T_2$, in the way
depicted in Fig.~\ref{fig-folding}. The construction continues
recursively, as values for $l_i$, $\gamma_i$ (and therefore $h_{i+1}
:= h_i / \cos \gamma_i$) are generated with each new trapezoid $T_i$.
We call the resulting region a \emph{polygonal spiral}, or simply
\emph{spiral}.

There are two of them, and they need not be equal, so we denote $N^R,
h_i^R, l_i^R, \gamma_i^R$, and $N^L, h_i^L, l_i^L, \gamma_i^L$, the
parameters of the right and the left spiral, respectively. These will
be determined later depending on $h_o$ and $l_o$, thus ultimately on
$k_f$. We will see to it that the following conditions hold:

\begin{itemize}
\item The spirals turn counterclockwise at each corner.

\item They have no self-intersections, or intersections between them
  or with the bulk of $\ta$.

\item For $\epsilon \in \{ R,L \}$, all angles $\gamma_i^\epsilon$ are
  rational multples of $\pi$.

\item There exists an absolute constant $K_1$ (i.e., $K_1$ does not
  depend on anything, including $k_f$) such that, for $\epsilon \in 
  \{ R,L \}$,
  \begin{equation}
    \label{cond-s1}
    h_{N^\epsilon}^\epsilon \le K_1 h_o.
  \end{equation}

\item There exists an absolute constant $K_2 > 1$ such that
  \begin{equation}
    \label{cond-s2}
    l_o \le \sum_{i=1}^{N^\epsilon} l_i^\epsilon \le K_2 \, l_o.
  \end{equation}

\item There exists an absolute constant $K_3$ such that, $\forall i=1,
  2, \ldots, N^\epsilon$,
  \begin{equation}
    \label{cond-s3}
    \frac{\tan \gamma_i^\epsilon} {l_i^\epsilon} \le \frac{K_3} 
    {h_i^\epsilon}.
  \end{equation}
  (The l.h.s.\ above is a measure of the ``curvature'' of the spiral 
  at the $i$-th corner.)
\end{itemize}

Under the above conditions the area of each spiral is bounded, as
$k_f \to 0$, because, dropping the superscript $\epsilon$,
\begin{eqnarray}
  \frac12 \sum_{i=1}^N ( 2l_i + h_i \tan \gamma_i ) h_i &\le& \frac{2 +
  K_3}2 \sum_{i=1}^N l_i h_i \nonumber \\
  &\le& \mbox{const } l_o h_N \\
  &\le& \mbox{const } l_o h_o \sim 1; \nonumber
\end{eqnarray}
having used, in this order, (\ref{cond-s3}), (\ref{cond-s2}), and
(\ref{cond-s1}).  Also, defining $(\ps, \ma, \mu)$ as in Section
\ref{sec-hyp}, namely, as the \dsy\ corresponding to the cross section
$\ps$ of all line elements based in $\bar{\bo} = \bo_f \cup \bo_d$, we
have:

\begin{proposition}
  \label{prop-hyp-s}
  $\ps$ is a global cross section for the \bi\ flow and $(\ps, \ma,
  \mu)$ is \hyp.
\end{proposition}

\proof First of all, $\ma$, as the first-return map onto $\ps$, is
well-defined almost everywhere (e.g., by the Poincar\'e Recurrence
Theorem). 

To prove that $\ps$ is a global cross section, we need to show that
a.a.\ \bi\ \tr ies have collisions against $\bar{\bo} = \bo_f \cup
\bo_d$. This is easy if we use a well-known result from the theory of
polygonal \bi s \cite{zk, bkm}: Let $P$ be the union of the two
spirals plus $R$, which is the rectangle (of base 1 and height $h$)
joining the open ends of the spirals. $P$ is a \emph{rational
polygon}, meaning that all its angles are rational multiples of
$\pi$. In a rational polygonal \bi, all but countably many values of
the velocity $u \in S^1$ are \emph{minimal}, in the sense that any
nonsingular flow-\tr y in configuration space (i.e., the set $\{ q(t)
\}_{t \in \R}$, provided that it contains no corner of $P$), with
initial velocity $u$, is dense in $P$ \cite{zk, bkm}.  This implies
that for a.a.\ initial conditions $(q,u)$, with $q \in P$, the \bi\
\tr y in $P$ hits the boundary of $R$, which means that the true
\bi\ \tr y, relative to the table $\ta$, hits $\bar{\bo}$.

As for the second assertion of Proposition \ref{prop-hyp-s}, we need
the following lemma, which will be proved later.

\begin{lemma}
  \label{lem-hyp-s}
  A material point that enters a spiral will travel all the way to the
  end of the spiral. In particular, if $\tau$ is the travel time
  between the last collision before entering the spiral and the first
  collision after exiting it (a.a.\ \tr ies eventually exit the
  spiral), then $\tau > 2l_o = 2 / k_f$.
\end{lemma}

Lemma \ref{lem-hyp-s} shows that Theorem \ref{thm-hyp} (and thus
Theorem \ref{thm-conhyp}) applies to the present case as well, since
its proof only requires of \tr ies visiting a strip---or a
spiral---that the travel time $\tau$ be larger than $2 / k_f$.  (Note
that, since the spirals are two polygons, they will have no effect on
the \hyp\ features of an infinitesimal beam of \tr ies, just like the
two strips. The only, inconsequential, difference is that the spirals
have more corners than the strips, resulting in more discontinuity
lines in $\ps$.)  
\qed

\proofof{Lemma \ref{lem-hyp-s}} The first assertion is an easy
consequence of our design, since a point that enters $T_i$ through the
shorter leg will necessarily exit it through the longer leg, thus
entering $T_{i+1}$ through the shorter leg, and so on. As for the
second assertion, clearly $\tau$ will be larger than twice the sum of
the lengths of the shorter bases of the trapezoids. By
(\ref{cond-s2}), this sum is bounded below by $l_o$.  
\qed

\fig{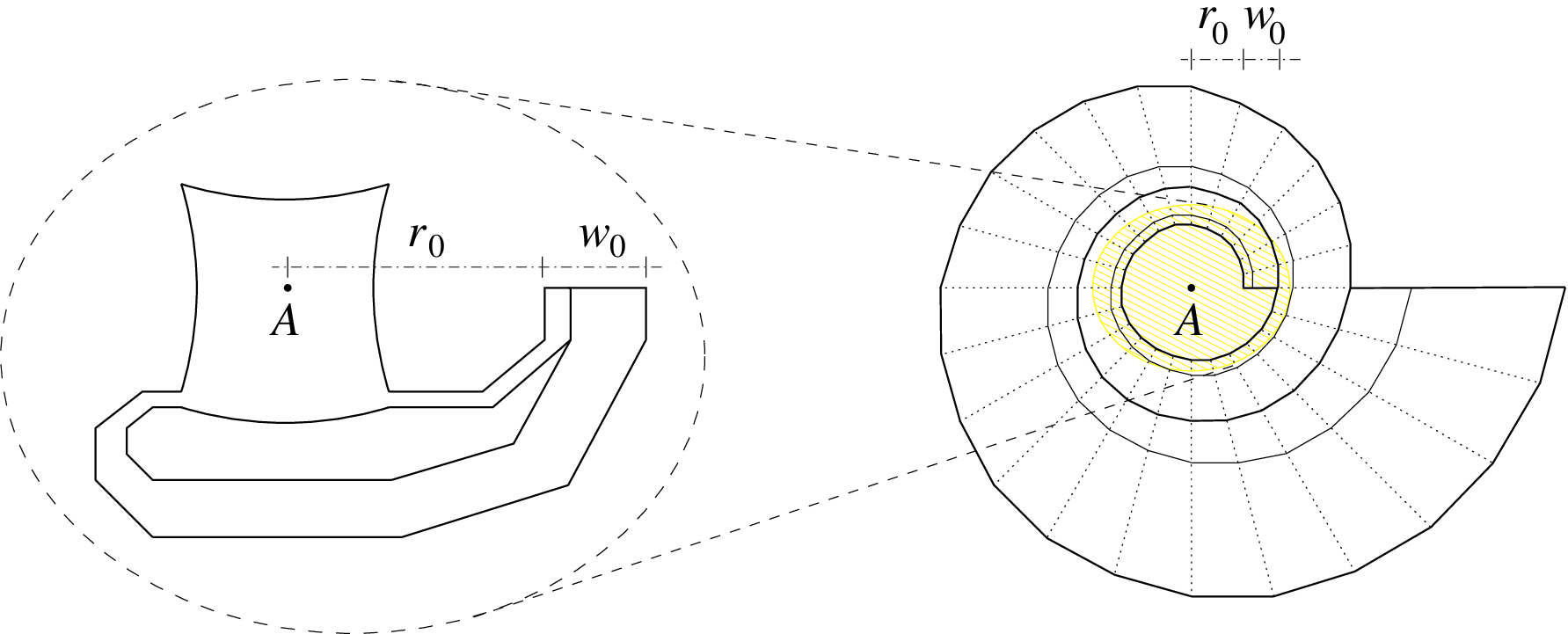} {13cm} {The double spiral (right picture) ``wrapping''
  around the bulk of $\ta$ (left picture). The double spiral starts
  when the two spirals coming out of the bulk of $\ta$ join. Its
  initial ray is $r_0$, its initial (total) width is $w_0$, each turn
  amounts to an angle $\bar{\gamma} = 2\pi / \bar{N}$, and the number
  of rounds is $M$. The point $A$ is the center of the double spiral.}

Let us finally give the exact construction of the two spirals. First
of all, we design the spirals to become adjacent after a finite number
of turns, say $m^R$ turns for the right spiral and $m^L$ turns for the
left spiral (left picture of Fig.~\ref{fig-reel}); $m^R$ and $m^L$ are
absolute constants.  We say that the two spirals have now joined in a
\emph{regular double spiral}, since they will keep adjacent as they
spiral outwards in the regular way shown in the right picture of
Fig.~\ref{fig-reel}.  More precisely, all trapezoids $T_i^R$, with $i
\ge m^R$, and $T_i^L$, with $i \ge m^L$, are similar, and are defined
by $\gamma_i^\epsilon = \bar{\gamma} := 2\pi / \bar{N}$, where
$\bar{N}$ is an integer (depending on $h_o$) to be determined
momentarily. The double spiral is also defined so that its initial ray
(meaning the distance from the border of the spiral to its center $A$,
see Fig.~\ref{fig-reel}) is $r_0$, an absolute constant so large that
intersection with the bulk of $\ta$ is avoided.

At each next corner, the ray (that is, the distance between that
corner and $A$) increases by a factor $1 / \cos \bar{\gamma}$.
Therefore, after the first round, the ray has become $r_{\bar{N}} :=
r_0 ( \cos \bar{\gamma} )^{-\bar{N}}$. Since the spiral wraps around
itself tightly (i.e., leaving no area uncovered), its initial width is
\begin{equation}
  \label{ds-10}
  w_0 := r_0 \left( \left(\cos \frac{2\pi} {\bar{N}}
  \right)^{-\bar{N}} \!\!\! - 1 \right).
\end{equation}
On the other hand, in the place where the left and right spirals join
to start the double spiral, one sees that
\begin{eqnarray}
  \label{ds-20}
  w_0 &=& h_{m^R}^R + h_{m^L}^L \nonumber \\
  &=& \left( \prod_{i=1}^{m^R} \frac1 { \cos \gamma_i^R } +
  \prod_{i=1}^{m^L} \frac1 { \cos \gamma_i^L } \right) h \\
  &=:& K_4 \, h . \nonumber
\end{eqnarray}
$K_4$ is an absolute constant if we prescribe that, for $i = 1,
\ldots, m^\epsilon$, the angles $\gamma_i^\epsilon$ are rational
multiples of $\pi$ and stay fixed while $k_f \to 0$ (this is
geometrically possible, cf.\ Fig.~\ref{fig-reel}, left picture). The
last two equations imply that
\begin{equation}
  \label{ds-30}
  h = h_1 = \frac{r_0} {K_4} \left( \left(\cos \frac{2\pi} {\bar{N}}
  \right)^{-\bar{N}} \!\!\! - 1 \right).
\end{equation}
Given $k_f$ sufficiently small, we use (\ref{ds-30}) to define both
$h$ and $\bar{N}$, keeping in mind that we want $h_o \le h \le K_1
h_o$, cf.\ (\ref{cond-s1}). We need this estimate from elementary
calculus:
\begin{equation}
  \label{ds-40}
  \lim_{n \to +\infty} \, \frac{n} {2\pi^2} \left( \left( \cos
  \frac{2\pi} {n} \right)^{-n} \!\!\! - 1 \right) = 1.
\end{equation}
So the r.h.s.\ of (\ref{ds-30}) decreases like $1/ \bar{N}$, as
$\bar{N} \to \infty$. This ensures that, given any sufficiently small
$h_o$, there exists an $\bar{N} = \bar{N}(h_o)$ such that the
corresponding $h = h(h_o)$, as in (\ref{ds-30}), verifies $h_o \le
h \le 2 h_o$.  Since $h_o = h_o(k_f)$, we rename these two values,
respectively, $\bar{N}(k_f)$ and $h(k_f)$ (abbreviated in $\bar{N}$
and $h$ when there is no risk of confusion). Clearly, as $k_f \to 0$,
\begin{eqnarray}
  \label{ds-50}
  && h(k_f) \sim h_o \sim k_f ; \\
  \label{ds-52}
  && \bar{N}(k_f) \sim h^{-1} \sim k_f^{-1} . 
\end{eqnarray}

Together with $r_0$, $h_{m^R}^R$ (equivalently $h_{m^L}^L$) and
$\bar{\gamma}$ (equivalently $\bar{N}$), the fourth and last parameter
that completely determines the double spiral is $M$, which is defined
as the number of complete rounds the spiral makes. (Once $M$ is
determined, the total number of trapezoids in the right and left
spirals is given by
\begin{equation}
  \label{ds-60}
  N^\epsilon = m^\epsilon + M \bar{N},
\end{equation}
for $\epsilon = R$ and $\epsilon = L$, respectively.) Choosing
\begin{equation}
  \label{ds-70}
  M = M(k_f) := \left[ \frac{l_o} {2\pi r_0} \right] + 1 = \left[
  \frac1 {2\pi r_0 k_f} \right] + 1
\end{equation}
(where $[ \,\cdot\, ]$ is the integer part of a positive number)
ensures that the first inequality of (\ref{cond-s2}) is verified, since
$\sum_i l_i^\epsilon > M 2\pi r_0 > l_0$. Also, for $\epsilon \in \{
R,L \}$,
\begin{equation}
  \label{ds-80}
  h_{N^\epsilon} = h_{m^\epsilon}^\epsilon (\cos 
  \bar{\gamma})^{-M \bar{N}} \sim h_{m^\epsilon}^\epsilon \sim h_o,
\end{equation}
as $k_f \to 0$, because of (\ref{ds-40}) and the fact that $M \sim
k_f^{-1}$ (whence $M \bar{N} \sim \bar{N}^2$). The above verifies
(\ref{cond-s1}). As for the second inequality of (\ref{cond-s2}), we
know that the trapezoids $T_i^\epsilon$, for $i \ge m^\epsilon$, are
similar.  Therefore, in the limit $k_f \to 0$, we obtain
\begin{eqnarray}
  \sum_{i=1}^{N^\epsilon} l_i^\epsilon &\sim& \sum_{i=m^\epsilon}^{N^\epsilon} 
  l_i^\epsilon \: = \ l_{m^\epsilon}^\epsilon \sum_{j=0}^{M \bar{N} - 1} 
  (\cos \bar{\gamma})^{-j} \nonumber \\
  &\sim& \tan \bar{\gamma} \, \frac{ (\cos \bar{\gamma})^{-M \bar{N}} 
  -1 } { (\cos \bar{\gamma})^{-1} -1 } \nonumber \\
  &\sim& \bar{N}^{-1} \frac1 { \bar{N}^{-2} } \sim \bar{N} \sim
  k_f^{-1} \\
  &\sim& l_o, \nonumber 
\end{eqnarray}
which proves (\ref{cond-s2}). In the above we have used (\ref{ds-52})
and the evident geometric equalities $l_{m^R}^R = r_0 \tan
\bar{\gamma}$ and $l_{m^L}^L = ( r_0 + h_{m^R}^R ) \tan \bar{\gamma}$
(Fig.~\ref{fig-reel}). Finally, (\ref{cond-s3}) holds because, for all
$i \ge m^\epsilon$, $l_i^\epsilon / h_i^\epsilon$ is constant, while
$\gamma_i^\epsilon = \bar{\gamma} \to 0$, as $k_f \to 0$.

The next and last result, whose proof is apparent, emphasizes the
motivation behind the constructions of Section \ref{sec-conf}.

\begin{proposition}
  The table $\ta = \ta(k_f)$ defined before is contained in a bounded
  region of the plane independent of $k_f$.
\end{proposition}

\end{document}